\documentclass[12pt]{amsart}
\usepackage{amssymb,amsmath,amsthm,amscd}
\usepackage[all]{xy}

\numberwithin{equation}{section}

\newtheoremstyle{fancy1}{10pt}{10pt}{\itshape}{12pt}{\textsc\bgroup}{.\egroup}{8pt}{
}
\newtheoremstyle{fancy2}{10pt}{10pt}{}{12pt}{\itshape}{.}{8pt}{ }

\theoremstyle{fancy1}

\newtheorem{cor}[equation]{Corollary}
\newtheorem{lem}[equation]{Lemma}

\newtheorem{thm}[equation]{Theorem}

\newtheorem{main}{Theorem}
\newtheorem*{main*}{Theorem}

\newtheorem*{cor*}{Corollary}

\theoremstyle{fancy2}

\newtheorem*{rem*}{Remark}

\newtheorem*{example*}{Example}

\newtheorem*{examples*}{Examples}

\newcommand{\cref}[1]{Corollary~\ref{#1}}

\newcommand{\lref}[1]{Lemma~\ref{#1}}

\newcommand{\tref}[1]{Theorem~\ref{#1}}

\newcommand{\gt}{\theta}

\newcommand{\Sph}{\mathbb{S}}

\newcommand{\C}{{\mathbb{C}}}
\newcommand{\R}{{\mathbb{R}}}
\newcommand{\Z}{{\mathbb{Z}}}

\renewcommand{\H}{\ensuremath{\operatorname{\mathsf{H}}}}

\newcommand{\G}{\ensuremath{\operatorname{\mathsf{G}}}}

\newcommand{\SO}{\ensuremath{\operatorname{\mathsf{SO}}}}
\renewcommand{\O}{\ensuremath{\operatorname{\mathsf{O}}}}

\newcommand{\SU}{\ensuremath{\operatorname{\mathsf{SU}}}}
\newcommand{\Spin}{\ensuremath{\operatorname{\mathsf{Spin}}}}

\newcommand{\T}{\ensuremath{\operatorname{\mathsf{T}}}}

\newcommand{\fg}{{\mathfrak{g}}}
\newcommand{\fk}{{\mathfrak{k}}}
\newcommand{\fh}{{\mathfrak{h}}}
\newcommand{\fm}{{\mathfrak{m}}}
\newcommand{\fn}{{\mathfrak{n}}}

\newcommand{\fp}{{\mathfrak{p}}}

\newcommand{\fso}{{\mathfrak{so}}}

\newcommand{\fsn}{ {\mathfrak{so}(n) } }

\newcommand{\pro}[2]{\langle #1 , #2 \rangle}

\def\con#1=#2(#3){#1 \equiv #2 \bmod{#3}}

\newcommand{\tr}{\ensuremath{\operatorname{tr}}}

\newcommand{\Ad}{\ensuremath{\operatorname{Ad}}}

\renewcommand{\sec}{\ensuremath{\operatorname{sec}}}

\DeclareMathOperator{\Id}{Id}

\newcommand{\no}{\noindent}
\newcommand{\co}{{cohomogeneity }}
\newcommand{\com}{{cohomogeneity one manifold }}

\newcommand{\Kpm}{\mathsf{K}^{\scriptscriptstyle{\pm}}}
\newcommand{\Kp}{\mathsf{K}^{\scriptscriptstyle{+}}}
\newcommand{\Km}{\mathsf{K}^{\scriptscriptstyle{-}}}

\newcommand{\Disc}{\mathbb{D}}

\newcommand{\gH}{\mathsf{H}}

\newcommand{\K}{\mathsf{K}}

\newcommand{\gK}{\mathsf{K}}
\newcommand{\gG}{\mathsf{G}}

\newcommand{\h}{h }

\begin{document}

\title{Non-negative curvature obstructions in cohomogeneity one
 and the Kervaire spheres}

\dedicatory{ Dedicated to Eugenio Calabi on his 80th birthday }

\author{Karsten Grove}
\address{University of Maryland\\
     College Park , MD 20742}

\email{kng@math.umd.edu}
\author{Luigi Verdiani}

\address{Universit\'a di Firenze\\ Via S. Marta 3\\50139 Firenze, Italy}
\email{verdiani@dma.unifi.it}

\author{Burkhard Wilking}
\address{University of M\"{u}nster\\
     Einsteinstrasse 62\\
     48149 M\"{u}nster, Germany}

\email{wilking@math.uni-muenster.de}
\author{Wolfgang Ziller}

\address{University of Pennsylvania\\
     Philadelphia, PA 19104}

\email{wziller@math.upenn.edu}
\thanks{The first named author was supported in part by the Danish Research
Council, the second named author was supported in part by a GNSAGA and the last
named author by the Francis J. Carey Term
Chair. The first and the last two authors were also supported by
grants from the National Science Foundation.}

\maketitle


Most examples of manifolds with non-negative sectional curvature are
obtained via product and quotient constructions, starting from
compact Lie groups with bi-invariant metrics. In \cite{GZ1} a new
large class of non-negatively curved compact manifolds was
constructed by using Lie group actions whose quotients are
one-dimensional, so called \co one actions. It was shown that if the
action has two singular orbits of codimension two, then the manifold
carries an invariant metric with non-negative sectional curvature.
This in particular gives rise to non-negatively curved metrics on all sphere
bundles over $\Sph^4$, which include the Milnor exotic spheres. In \cite{GZ1} it was also conjectured
that any \co one manifold should carry an invariant metric with
non-negative sectional curvature. The most interesting application
of this conjecture would be the Kervaire spheres, which can be presented as particular   $2n-1$
dimensional Brieskorn varieties $M^{2n-1}_d \subset \mathbb{C}^{n+1}$
 defined  by the equations
$$ z_0^d + z_1^2 +\cdots z_n^2=0 \quad , \quad |z_0|^2 + \cdots |z_n|^2 =1 .$$

For $n$ odd and $d$ odd, they are homeomorphic to spheres, and are
exotic spheres if $\con 2n-1=1(8)$. As was discovered by E.Calabi in
dimension 5 and later generalized in \cite{HH}, the Brieskorn
variety $M^{2n-1}_d$ carries
 a \co one action by \SO(2)\SO(n)
defined by $(e^{i\gt},A)(z_0,\cdots ,
z_n)=(e^{2i\gt}z_0,e^{id\gt}A(z_1,\cdots,z_n)^t)$. It was observed
in \cite{BH} (see \cite{searle} for $n=3$) that $M^{2n-1}_d$ does
not admit an $\SO(2)\SO(n)$ invariant metric with positive
curvature, if $n\ge 3 \, , \,d\ge 2$. On the other hand, the
non-principal orbits of this action have codimension 2 and $n-1$.
Hence,  if $n=3$, they carry an invariant metric with non-negative
curvature. If $n=2$ the action is a linear action on the lens space
$L(1,d)$. If $d=1$ the action is the linear tensor product action on
$\Sph^{2n-1}$ and if $d=2$ it is a linear action on $\T_1 \Sph^n$.
Hence in all these cases one has an invariant metric with
non-negative curvature. In all remaining cases we prove:

\begin{main}
For $n\ge4$ and $d\ge 3$, the Brieskorn variety $M^{2n-1}_d$ does not support an $\SO(2)\SO(n)$ invariant metric with non-negative sectional curvature.
\end{main}

If $n=4$, the normal $\SU(2)$ subgroup in $\SO(4)$ acts freely on
$M^7_d$. Its quotient is $\Sph^4$, since the induced \co action on
the base is the sum action of $\SO(2)\SO(3)$ on $\Sph^4$. In
\cite{GZ1} it was shown that such principal bundles admit another
cohomogeneity one action by $\SO(4)$ with non-principal orbits of
codimension two and hence an invariant metric with non-negative
sectional curvature. We do not know if any of the remaining
manifolds $M^{2n-1}_d$ admit non-negative curvature.

\bigskip

In addition to the cohomogeneity one Brieskorn varieties,  we
construct several other infinite families of cohomogeneity one
manifolds that cannot support invariant metrics of non-negative
curvature (see Theorem \ref{newones}). In particular, we get such an
infinite family for any codimensions bigger than two   (see
Corollary \ref{maincor}). Together with the examples of the
Brieskorn varieties
 we obtain the following counterpart
to the positive result of \cite{GZ1}:

\begin{main}
For any pair of integers $(\ell_-,\ell_+) \ne (2,2)$ with
$\ell_{\pm} \ge 2$, there is an infinite family of cohomogeneity one
manifolds $M$ with singular orbits of codimensions $\ell_{\pm}$ and
no invariant metric of non-negative sectional curvature.
\end{main}

\bigskip

Our results are in contrast to the homogeneous case, where all
compact manifolds admit an invariant metric of non-negative
curvature, and they leave the classification of non-negatively
curved cohomogeneity one manifolds wide open. We point out though
that all cohomogeneity one manifolds admit invariant
 metrics with almost non-negative sectional curvature (see \cite{ScTu}).
 We also point out that our obstructions are
global in nature since each half, as homogeneous disc bundles over
the singular orbits, naturally admit invariant metrics with
non-negative curvature.

\bigskip

In Section 1 we describe the \co one action on $M^{2n-1}_d$ and the
invariant metrics. In Section 2 we derive the obstructions to
non-negative curvature and prove Theorem A.  In section 3 we exhibit
other infinite families of \co manifolds that do not admit invariant
metrics with non-negative curvature.

\section{Cohomogeneity one Brieskorn varieties and invariant metrics}

For any integer $d\geq 1$ the Brieskorn variety $M_d^{2n-1}$ is
the smooth $(2n-1)$-dimensional submanifold of $\mathbb
C^{n+1}$ defined by the equations

$$ \left\{
\begin{array}{l}

z_0^d+z_1^2+z_2^2+\ldots+z_n^2=0,\\
|z_0|^2+|z_1|^2+\ldots+|z_n|^2=1.
\end{array}\right.$$

When $d=1$,  $M_1^{2n-1}$ is diffeomorphic to the standard sphere
via $(z_0,Z)\to Z/|Z|$. When $d=2$, $M_2^{2n-1}$ is diffeomorphic to
the unit tangent bundle of $S^{n}$ since $X+iY$ lies in $M_2^{2n-1}$
iff $\pro{X}{X} = \pro{Y}{Y} =  \frac{1}{2} , \pro{X}{Y} = 0$. When
$n=2$, $M_d^3$ is the Lens space $L(1,d)$  while, for $n\geq 3$ the
Brieskorn varieties are simply connected.

\smallskip

 The group $\G=\SO(2) \SO(n)$ acts on $M_d^{2n-1}$ by
$$(e^{i\theta},A) (z_0,Z)=(e^{2i \theta} z_0, e^{i d \theta} A Z)
\, , \, (z_0,Z)\in \C\oplus\C^n.$$
 $|z_0|$ is invariant under this action
and one easily sees that two points  belong to the same $\G$-orbit
if and only if they have the same value of $|z_0|$. Furthermore,
using the Schwarz inequality, we have that $0\leq |z_0|\leq t_0$,
where $t_0$ is the unique positive solution of $t_0^d+t_0^2=1$. In
particular, the manifold is \co one and  $M/\G=[0,t_0]$.

\smallskip

In order to describe the possible orbit types for this action, we
consider the Lens space $M_d^3\subset M_d^{2n-1}$, defined by
$$M_d^3=\{(z_0,z_1,z_2,0,\ldots,0)\}\cap M_d^{2n-1}.$$
The points of $M_d^3$ such that $0< |z_0|< t_0$ have a common
isotropy subgroup, namely
$$\H=\Z_2\times \SO(n-2)=\begin{cases}

(-\epsilon, \text{diag}(\epsilon,\epsilon,A)) &\text{if $d$ odd}\\
(\epsilon, \text{diag}(1,1,A)) &\text{if $d$ even}

\end{cases}
$$
where $\epsilon=\pm 1$ and $A\in \SO(n-2)$. Conversely, $M_d^3$
coincides with the fixed point set of $\H$. There are three orbit
types, the principal ones, corresponding to the interior of the
interval $[0,t_0]$, and two singular ones corresponding to the
boundary points.

One singular orbit,  denoted  by $B_-$,  is  given by $z_0=0$ and
has codimension $2$. Choosing $p_-=(0,1,i,0,\ldots,0)\in M_d^3$,
the isotropy subgroup of $p_-$ is
$$\Km=\SO(2)\SO(n-2)=(e^{-i\theta},\text{diag}(R(d\theta),A)).$$
where $R(\theta)$ is a counterclockwise rotation with angle
$\theta$.

The second singular orbit, $B_+$,  is given by $|z_0|=t_0$ and has
codimension $n-1$. If we choose
$p_+=(t_0,i\sqrt{t_0^d},0,\ldots,0)$, the corresponding isotropy
subgroup  is
$$\Kp=\begin{cases}

\O(n-1)=(\det B, \text{diag}(\det B,B)) &\text{if $d$ odd}\\
\mathbb Z_2\times \SO(n-1)=(\epsilon, \text{diag}(1,B')) &\text{if
$d$ even}

\end{cases}
$$
where $\epsilon=\pm 1$, $B\in \O(n-1)$ and $B'\in \SO(n-1)$.

\smallskip

A tubular neighborhood of the singular orbit $\G/\Kpm$  is
equivariantly diffeomorphic to the homogeneous vector bundle
$\G\times_{\Kpm} \Disc_\pm$, where $\Disc_\pm$ can be identified
with the normal disc to the orbit $B_\pm=\G p_\pm = \G/\Kpm $ at
$p_\pm$. The group $\Kpm$ acts linearly on $\Disc_\pm$,  and
transitively on the sphere $\partial \Disc_{\pm}=\Sph_{\pm}=\Kpm/\H
$. Furthermore, $M$ is the union of these tubular neighborhoods  $M
= \G \times_{\Km} \Disc_- \cup \G \times_{\Kp} \Disc_+$ glued along
a principal orbit $\G/\H$. The isotropy groups $\Kpm$ are isotropy
groups of the end points of a minimal geodesic joining the singular
orbits for an invariant metric on $M$.

\smallskip

The largest normal subgroup common to $\G$ and $\Km$, which is the
ineffective kernel of the $\G$ action on $B_-$, is the cyclic group
generated by $(e^{-i\pi/d},-\Id)$ for n even and
$(e^{-2i\pi/d},\Id)$ for n odd. On the other hand, the action of
$e^{i\theta}\in\SO(2)\subset\Km$ on the slice $\Disc_-$ is given by
$R(2\theta)$ since $e^{i\pi}\in\H$.
 Hence the singular orbit $B_-$ is the fixed point set of
 a group of isometries and thus is totally geodesic if $d\ge 3$.
  This implies, as was
observed in \cite{BH}, that $M$ cannot carry a $\G$ invariant metric
with positive curvature if $n\ge 4$ and $d\ge 3$ since $B_-$ does
not admit a homogeneous metric with positive curvature.

\bigskip

We now  describe the Riemannian metrics $g$ on $M_d^{2n-1}$
invariant under the action of  $\G=\SO(2)\SO(n)$. Let $\gamma$ be a
geodesic orthogonal to one and hence all $\G$ orbits, which we can
assume lies in $M^{\H}$ and ends at $p_+$ chosen as above. Its
beginning point  $\gamma(0)$ will lie in $B_-\cap M^{\H}$, which
consists of the two circles $\{ (0,z,iz)\mid |z|=1\}$ and $\{
(0,z,-iz)\mid |z|=1\}$. By changing the sign of the third basis
vector if necessary, we can assume that $\gamma(0)$ lies in the
first circle where the isotropy group is constant equal to $\Km$. We
can hence assume, changing $p_-$ if necessary, that $\gamma$ is
parameterized by arc length with $\gamma(0)=p_-$ and $\gamma(L)=p_+$
and all isotropy groups are described as above.

For $0<t<L$,
 $\gamma(t)$ is a regular
 point  with constant isotropy group $\H$ and the metric on the
 principal
 orbits $\G \gamma(t)=\G/\H$ is a family of homogeneous metrics
 $g_t$.
 Thus on the regular part
 the metric is determined by
$$g_{\gamma(t)}=d\,t^2+g_t$$
and since the regular points are dense it also describes the metric
on $M$.

By means of  Killing vector fields, we identify the tangent space to
$\G/\H$ at $\gamma(t), t\in(0,L)$ with an $\Ad(\H)$-invariant
complement $\fn$ of the isotropy subalgebra $\fh$ of $\H$ in $\fg$
and the metric $g_t$ is identified with an $\Ad(\H)$ invariant inner
product on $\fn$. The complement $\fn$ can be chosen as the
orthogonal complement to $\fh$ with the respect to a fixed
$\Ad(\H)$-invariant scalar product $Q$ in $\fg$. In our case a convenient
choice of $Q$ on $\mathfrak{so}(2)+\mathfrak{so}(n)$ is
$$Q(a+A,b+B)=d^2ab -\frac{1}{2} \tr(A\cdot B)$$
where $a, b\in\fso(2) $ and  $A,B\in \fsn$. We obtain a
$Q$-orthogonal decomposition:
$$T_{\gamma(t)} \G(\gamma(t))\simeq \fn \simeq \fp_1+\fp_2+\fm_1+\fm_2$$
for $t\in(0,L)$, where $\fp_i$ are one-dimensional trivial
$\Ad(\H)$-modules and $\fm_i$ are   $(n-2)$-dimensional irreducible
$\Ad(\H)$-modules. A $Q$-orthonormal basis of $\fn$ is given by
$$X=(\frac 1d I+E_{12})/\sqrt{2}\in\fp_1 \, ,\, Y=(\frac 1d
I-E_{12})/\sqrt{2}\in\fp_2\, ,\, E_i=E_{1(i+2)}\in\fm_1\, ,\,
F_i=E_{2(i+2)}\in\fm_2$$   for $ i=1..n-2$, where $E_{ij}$ is the
standard basis for $\mathfrak{so}(n)$ and $I=E_{12}\in\fso(2) $.

The same identification at the singular points takes the form
$$\fk_-=\fh+\fp_1,\qquad T_{\gamma(0)} B_-\simeq \fp_2+\fm_1+\fm_2,$$
$$\fk_+=\fh+\fm_2,\qquad T_{\gamma(L)} B_+\simeq \fp_1+\fp_2+\fm_1.$$

Any $\Ad(\H)$-invariant scalar product on $\fn$ must be a multiple of
$Q$ on the irreducible $\H$-modules $\fm_i$ and the inner products
between $\fm_1$ and $\fm_2$ are described in terms of an
$\Ad(\H)$-equivariant isomorphism. Up to a factor, the only such
isomorphism $\fm_1\to\fm_2$ is given by $E_i\to F_i$. Furthermore,
$\fp_1 +\fp_2$ must be orthogonal to $\fm_1+\fm_2$. In terms of this
basis, $g_t$ can  therefore be described as follows:

$$g_t(X,X)=f_1^2(t),\qquad g_t(Y,Y)=f_2^2(t),\qquad g_t(X,Y)=f_{12}(t)$$

$$g_t(E_i,E_i)=\h_1^2(t),\qquad g_t(F_i,F_i)=\h_2^2(t),\qquad g_t(E_i,F_i)=
\h_{12}(t)$$ with all other inner products being 0.

\smallskip

If these 6  functions are smooth, and $g_t$ is positive definite,
they define a smooth metric on the
 regular part of $M$. Some further conditions are required to ensure
that the metric $g_t+d\,t^2$ can be smoothly extended to the singular
 orbits. The precise necessary and sufficient conditions were determined in
 \cite{BH}.
 For our purposes
it is sufficient to consider  only a few of these smoothness
conditions. They are, assuming that $d\ge 3$:

\begin{eqnarray}\label{smooth}
 f_1(0)=0 \,\, ,\,  f_1'(0)=\sqrt{2}/d  \,\, ,\,
 f_{12}(0)=f_{12}'(0)=0\qquad\qquad\qquad
 \\ \nonumber
\h_1(0)= \h_2(0)   \,\, ,\,  \h_1'(0)= \h_2'(0)=0 \,\, ,\,
\h_{12}(0)=\h_{12}'(0)=0 \,\, ,\, \h_1'(L)=0 \,\, ,\, \h_2(L)=0
\end{eqnarray}

 The first two follow from the fact that  $\fp_1$ collapses at
$B_-$ and that the circle $\Km/\H$ has length $2\pi d/\sqrt{2}$ in
the metric $Q$. The inner product $f_{12}(0)=0$ since the slice $\Disc_-$ and $B_-$
are orthogonal.  Since $\Km$ fixes $Y$, it follows that
  $\langle X,Y\rangle$, restricted to the slice, is a function of
  the square of the distance to the origin. Hence $f_{12}(t)$ is an even
function of $t$ which implies that
 $f_{12}'(0)=0$. We have $\h_{12}(0)=0$ and $\h_1(0)= \h_2(0)$ since the isotropy representation
of $\Km$ on $\fm_1+\fm_2$ is irreducible. The remaining derivatives
vanish at $t=0$ since $B_-$ is totally geodesic. Now  $\h_1'(L)=0$ means
that the second fundamental form of $B_+$ restricted to $\fm_1$ is
0. Indeed, this follows from the fact that the second fundamental
form is equivariant under the action of $\Kp$, and that this action
on the slice $\Disc_+$ is equivalent to its action on the subspace
of $T_{p_+}(B_+)$ spanned by $X+Y=E_{12}$ and $\fm_1$. Finally,
$\h_2(L)=0$ since $\fm_2$ collapses at $p_+$.

\section{Obstructions and Proof of Theorem A}

From now on we assume that $M=M_d^{2n-1}$ with $d\ge 3$ and $n\ge 4$
and that $M$ is equipped with a $\G=\SO(2) \SO(n)$-invariant
$C^2$-metric $g$. Our starting point is the following global fact:

\begin{lem}\label{constant}
If $(M,g)$ is non-negatively curved, then $\h_{12}(t)=0$  and
$\h_1(t)$ is constant for $t\in[0,L]$.
\end{lem}

\begin{proof}

Let $E_i\in \fm_1$ as before and  consider the vector field $\tilde{E}_i$
obtained by parallel transport of $E_i(\gamma(L))$ along $\gamma(t)$
for $t\in [0,L]$. By $\Ad(\H)$-equivariance,   the parallel transport
along $\gamma$ preserves $\fm_1+\fm_2$ and hence $\tilde
E_i(\gamma(0))\in \fm_1+\fm_2$. Recall that $B_-$ is totally
geodesic and that the second fundamental form of $B_+$ restricted to
$\fm_1$ is 0. Since the normal geodesic connecting $p_-$ to $p_+$ is
a minimizing connection between $B_-$ and $B_+$, the second
variation formula implies that $\pro{R(\tilde E_i,
\gamma')\gamma'}{\tilde E_i} = 0$. Since $g$ has non-negative
curvature, $R(\tilde E_i, \gamma')\gamma'=0$ and hence $\tilde E_i$
is a Jacobi field along $\gamma$, which has the same initial
conditions as $E_i$ at $\gamma(L)$. Thus both must agree, which
means that $E_i$ is parallel along $\gamma$ and hence $\h_1$ is
constant for $t\in [0,L]$.

The orthogonal complement $\fm_1^\perp$ of $\fm_1$ in $\fm_1+\fm_2$
is both invariant under parallel translation and under the symmetric
endomorphism $R(\cdot , \gamma')\gamma'$. Therefore any Jacobi field
with initial conditions in $\fm_1^\perp$ stays in $\fm_1^\perp$. By
\eqref{smooth}   we have $\fm_1^\perp=\fm_2$ for $t=0$ and hence
$\fm_1$ and $\fm_2$ are orthogonal everywhere.
\end{proof}

 The
remainder will be a proof by contradiction. By computing the
sectional curvature of two particular 2-planes, the non-negativity
condition gives   upper and lower bounds for $\h_2'(t)$ for $t$
small, which when combined yields $d\le 2$.

\smallskip

We can replace all functions, e.g. $\h_1(t)$, by $a\h_1(t/a)$
corresponding to multiplying the metric by $a$.  From now on we will
normalize the metric so that $\h_1(0)=1$.  It follows from
\eqref{smooth} that $\h_2(0)=1$ as well.

\smallskip

 Let $A,B\in\fn$ be linearly independent.  The corresponding
  Killing vector fields  span a two
 plane tangent to the principal orbits  at the
 regular points of $\gamma(t)$.
  The sectional curvature of such a two plane can be
   computed
  using the Gauss equations
 $$R(A,B,B,A)=\overline R(A,B,B,A)-\tfrac{1}{4}\, {g_t(A,B)'}^2+\tfrac{1}{4}\,
 g_t(A,A)'\,g_t(B,B)'$$
 where $\overline R(A,B,A,B)$ is the intrinsic curvature of the principal orbit
  $\G/\H$. The latter can be computed by any of the well known
 formulas for  the curvature tensor of a homogeneous space
 (cf. \cite[7.30]{Be} or \cite{GZ2}).

\smallskip

 We choose the vectors $A=E_1+F_2$
  and $B=E_2+F_1$, which is possible since $n\ge 4$.
  Using the fact that $[A,B]=0$, a curvature computation
  shows:

 \begin{equation*}\label{non-neg}
 R(A,B,B,A)= \tfrac{1}{2}\,(1-\h_2^2)^2
 \frac{f_1^2+f_2^2-2\, f_{12}^2} {f_1^2\,f_2^2-f_{12}^2}-\h_2^2\, {\h_2'}^2
 \end{equation*}

  \no
In the following we set

 \[
\delta(t)=1-\h_2(t)^2
\]

\no and hence the non-negative curvature condition can be formulated
as:

\begin{equation}\label{cond1mod}
 (\delta')^2 \le
 2\frac {f_1^2+f_2^2-2\,f_{12}^2}{(f_1^2\, f_2^2- f_{12}^2)}
 \delta^2
 \end{equation}

Since  $F_1/\h_2$ is a parallel vectorfield, it follows that
 $\sec(\gamma'(t),F_1)=-\h_2''/\h_2$, and hence the function $\h_2(t)$ is
 concave. Using the smoothness condition $\h_2'(0)=0$ and $\h_2(0)=1$
  it follows that
   $\delta(t)\ge 0$. We first claim that $\delta(t)> 0$ for $t> 0$.
   If not, there would be a $t_0 > 0$ with $\delta(t_0)=0$ but
   $\delta(t)>0$ for $t>t_0$ since $\h_2(L)=0 $ (cf. \eqref{smooth}). But then
    \eqref{cond1mod} implies that
   $\delta'/\delta= \log(\delta)'$ is bounded from above, which by integrating from
   $t_0+\epsilon$ to $t_0+\eta$ leads to a
   contradiction when $\epsilon\to 0$.

The smoothness conditions \eqref{smooth} imply that $f_1(0)=0$ and
$\lim_{t\to 0}\tfrac{f_{12}}{f_1}=0$ and since $f_2(0)\ne 0$,
\eqref{cond1mod} implies that:
\begin{equation}\label{eq: first}
 \log(\delta)'\le (1+\eta)\frac {\sqrt{2}}{f_1}\mbox{ for $t\in (0,L)$}
 \end{equation}
 where $\eta(t)$ is a positive function with $\lim_{t\to 0}\eta(t)=0$.

\smallskip

To obtain a lower bound, we consider the family of $2$-planes
spanned by $A_r=X+r\,Y$ and
 $B=F_1$.
 A necessary condition for $R(A,B,B,A)$ to be non-negative for all $r$ is
 \begin{equation*}
 \label{cond3}
 R(X,F_1,X,F_1)\,R(Y,F_1,Y,F_1)\geq R(X,F_1,Y,F_1)^2
 \end{equation*}

\no A computation shows that
 \begin{eqnarray*}
 R(X,F_1,X,F_1)&=&
 -\tfrac{1}{8}\delta\bigl(4-2f_1^2-2f_{12}-\delta\bigr)+
 \tfrac{1}{8}(f_1^2+f_{12})^2
 +\tfrac{1}{2}f_1 f_1'\delta'\\
 R(Y,F_1,Y,F_1)&=&
 -\tfrac{1}{8}\delta\bigl(4-2f_2^2-2f_{12}-\delta\bigr)+
 \tfrac{1}{8}(f_2^2+f_{12})^2
 +\tfrac{1}{2}f_2 f_2'\delta'\\
 R(X,F_1,Y,F_1)&=& -\tfrac{1}{8}\delta
 \bigl(4-f_1^2-f_2^2-2f_{12}-\delta\bigr)
+\tfrac{1}{8}(f_1^2+f_{12})(f_2^2+f_{12})
+\tfrac{1}{4}f_{12}'\delta',
\end{eqnarray*}
 If we set $a=f_2(0)^2$, and using $\delta(0)=\delta'(0)=0$,
  it follows that
\begin{eqnarray*}
R(X,F_1,X,F_1)\,R(Y,F_1,Y,F_1)-R(X,F_1,Y,F_1)^2=
 -\tfrac{a^2}{16}(1+\eta_1)\delta+
 \tfrac{a^2}{16}(1+\eta_2)f_1f_1'\delta' \ge 0
 \end{eqnarray*}

\no  where $\eta_i$ are functions with $\lim_{t \to 0}\eta_i(t)=0$,
$i=1,2$. Combining with inequality~\eqref{eq: first} gives
\[
(1-\eta_3(t))\tfrac{1}{f_1(t)f_1'(t)}\leq \log(\delta)'\leq
(1+\eta_3(t))\frac {\sqrt{2}}{f_1(t)}
\]
for a suitable positive function $\eta_3$ with $\lim_{t \to
0}\eta_3(t)=0$. In particular $ f_1'(0)\ge \tfrac{1}{\sqrt{2}}$ and
since  $f_1'(0)=\tfrac{\sqrt{2}}{d} $ (cf.\eqref{smooth}), we obtain
$d\le 2$. This finishes the proof of our main Theorem.

\begin{rem*}
The  first inequality is already enough to rule out the existence of
invariant non-negatively curved metrics of class  $C^{d+1}$ for
$d=2m+1>2$. Indeed, since $f_1'(0)=\tfrac{\sqrt{2}}{d}$, \eqref{eq:
first} implies that
 \[
 \log(\delta)'\le (1+\tilde{\eta}(t))\frac {d}{t}.
\]
Integrating gives that  for each $\alpha> d$ we have $\delta(t)\ge
t^\alpha$
 provided that $t$ lies in a sufficiently small interval
 $(0,\varepsilon(\alpha))$.
On the other hand, using the action of the cohomogeneity one Weyl group at $B_-$ (see e.g. \cite{AA}), one
shows that if $d$ is odd $\delta(t)=0$ for $t\in [-L,0]$. Thus the
metric can not be of class $C^{k}$ for any integer $k>d$.
\end{rem*}

The choice of the given $2$-planes can be  motivated as
 follows. For the first set of 2-planes, the
  contribution given by the second fundamental form is
 non-positive. This forces the intrinsic sectional curvature to
 be non-negative.
Since $[A,B]=0$, the 2 plane has intrinsic curvature 0 in a
bi-invariant metric, and non-negative curvature of $M$ gives a strong
upper bound for $\h_2'$.

For the second family, we choose 2-planes whose intrinsic sectional
curvature should be negative, since it forces the contribution from
the second fundamental form to be positive, giving a lower bound for
$\h_2'$. This can be done by considering two planes containing
$E_{12}\in\fso(2)$ since non-negative curvatures for such 2-planes
would imply that $\SO(2)$ splits off isometrically, contradicting
the fact that
 $X=\frac 1d A+E_{12}$ must vanish at
 the singular orbit $B_-$.

 \section{Other Examples and Proof of Theorem B}

It  turns out that there are many more examples of compact
cohomogeneity one manifolds which do not admit invariant non-negatively
curved metrics. Recall that any choice of groups $\H \subset \{\Km, \Kp\} \subset \G$ with $\Kpm/\H = \Sph^{l_{\pm}}$ define a cohomogeneity one manifold $M = \G \times_{\Km} \Disc^{l_{-} + 1} \cup  \G \times_{\Kp} \Disc^{l_{+} + 1}$.

 Suppose $\gK'/\gH'=$$ \Sph^{\ell -1}$ is a sphere, and
let $sl\colon \gK'\rightarrow \O(\ell)$ be the corresponding
representation. Suppose furthermore
 that $\mu\colon \gK'\rightarrow \SO(k)$
is an irreducible faithful representation with $\mu(\gH')\subset
\SO(k-1)$. We can then define a cohomogeneity one manifold $M$ with
$\G:=\SO(n)$. View $\SO(k)$ as the upper $k\times k$-block  of
$\SO(n)$ and set

\begin{align}\label{ex}\nonumber
\Km &= \mu(\gK')\cdot \SO(n-k)\subset \SO(k)\SO(n-k)\subset \SO(n)\\
 \Kp &=\mu(\gH')\cdot \SO(n-k+1) \subset
\SO(k-1)\SO(n-k+1)\subset \SO(n)\\\nonumber \H\;\, & =
\mu(\gH')\cdot \SO(n-k)\subset\Kpm
\end{align}

\begin{thm}\label{newones}   Let $\gK'/\gH'=$$ \Sph^{\ell -1}$ with $\ell\ge 3$
and $\mu\colon \gK'\rightarrow \SO(k)$  an irreducible faithful
representation  such that

\begin{itemize}
\item[a)]
 $\mu(\gH')\subset
\SO(k-1)$
\item[b)]
$\mu(\gK')$ does not act transitively on $\Sph^{k-1}$
\item[c)]
$sl$ is not a subrepresentation of $S^2\mu$
\item[d)]
 $n\ge k+2$
\end{itemize}
Then $M$ as defined by \eqref{ex} does not admit a $\G$ invariant
$C^2$ metric of non-negative sectional curvature.
\end{thm}

\begin{proof} We assume that there is an invariant non-negatively
curved metric on $M$. Let $c\colon [0,L]\rightarrow M$ denote a
shortest geodesic from the singular orbit $\gG/\gK_-$ to $\gG/\gK_+$
and  extend $c$ to a geodesic $c\colon \R\rightarrow M$. We can
assume that $\H$ leaves $c$ pointwise fixed. Let $W\subset
\fn=\fh^\perp\subset\fso(n)$ be the subspace spanned by $E_{i,j}$
with $1\le i\le k \, ,\, k+1\le j\le n$ and $W(t)$ its image in
$T_{c(t)}M$. Then $W(t)$ is invariant under parallel translation since it
is complementary and hence orthogonal to the fixed point set  of
$\SO(n-k)\subset\H$ in $\fn$.

Notice that $W(0)$ and $W(2L)$ are tangent to $B_-=\G/\Km$ and we
first claim that the second fundamental form of $B_-$ vanishes on
them. Indeed, the action of $\Km$ on $W$ is given by
$\mu\hat{\otimes}\rho_{n-k}$, where as usual we denote by
$\rho_{n-k}$ the representation of $\SO(n-k)$ on $\R^{n-k}$. Hence
$S^2(W)= S^2(\mu)\hat{\otimes}S^2(\rho_{n-k}) \oplus
\Lambda^2(\mu)\hat{\otimes}\Lambda^2(\rho_{n-k}) $. Furthermore,
$\Lambda^2(\rho_{n-k})\simeq\fso(n-k)$ is irreducible and
$S^2(\rho_{n-k})$ is a trivial representation plus an irreducible
one. The claim follows from equivariance of the second fundamental
form  under the action of $\Km$ since the slice representation of
$\Km$ is given by $sl\hat{\oplus}\Id$ and by assumption $sl\notin
S^2(\mu)$.

We now consider the geodesic $c_{[0,2L]}$ as a point in the space of
all path in $M$ which start and end in $\gG/\gK_-$. Clearly $c$ is a
critical point of the energy functional and we claim that it has
index $n-k$ with respect to any $\G$-invariant metric. To see this,
observe that since geodesics from $\gG/\gK_-$ to $\gG/\gK_+$ are
minimal, this path space can be approximated by broken geodesics
starting and ending orthogonal to  $\gG/\gK_-$ and with break points
in $\gG/\gK_+$ without changing the index and nullity of the
critical point $c$. On this subspace of dimension
$\dim(\gK_-/\gH)+\dim(\gK_+/\gH)$  the energy clearly has nullity
$\dim(\gK_-/\gH)$ and index $\dim(\gK_+/\gH)=n-k$.

On the other hand, we can extend each vector in $W(0)$  to a
parallel vectorfield along $c$ and by the previous observation the
boundary terms in the variational formula for the energy vanish.
Hence  the Hessian of the energy is negative semidefinite on this
$k(n-k)$-dimensional vectorspace. Since the index is only $(n-k)$ it
follows that there exists a subspace $W_1\subset W$ of dimension
 at least $(k-1)(n-k)$ which is contained in the nullspace of the
 Hessian and hence consists of parallel Jacobifields.
 As in
\lref{constant}, it follows that these vectorfields are actually
given by Killing fields. Since $W_1$ is also invariant under the
action of $\SO(n-k)\subset\H$ we can assume, via a change of basis,
that $W_1$ is spanned by the first $k-1$ rows in $W$. Let $W_2$ be
last row in $W$. Since $\mu$ is irreducible, it follows that
$W_1(0)$ is orthogonal to $W_2(0)$ and hence, as in the proof of
\lref{constant}, $W_1(t)$ is orthogonal to $W_2(t)$ for all $t$.

Using again that $\mu$ is irreducible, we can normalize the metric
so that $E_{i,j}\in W$ is an orthonormal basis at $t=0$. By the
above it follows that $E_{i,j}\in W_1$ remains orthonormal for all
$t$. Since $\SO(n-k)\subset\H$ acts irreducibly on $W_2$, we can set
$\h(t):=\|E_{i,j}\|$ for $E_{i,j}\in W_2$. non-negative curvature
implies that the function $\h$ is concave, and since $\h(0)=1$,
$\h'(0)=0$, it follows that $\h(t)\le 1$. As before we set
$\delta=1-\h^2\ge 0$.

Let  $Y_1=E_{k,k+1}$ and $Y_2=E_{k,k+2}$ and  set $X_1=\Ad_a
E_{k-1,k+1}$ and $X_2= \Ad_a E_{k-1,k+2}$,  where $a\in\SO(k-1)$ is
an element in the upper $(k-1)\times (k-1)$--block. We now consider
the 2-plane spanned by  $A=X_1+Y_2$ and $B=X_2+Y_1$. Clearly $A$ and
$B$ commute. A computation shows that:

 \begin{equation}\label{Bound}
 R(A,B,B,A)= \frac 1 4 \delta^2 Q([Y_2, X_2] -  [X_1, Y_1],
 P^{-1}_t( [Y_2, X_2] -  [X_1, Y_1] ) )-\tfrac{1}{4}\delta'^2\ge 0
\end{equation}

\no where $Q(X,Y)=-\frac 1 2 \tr(XY)$ is a bi-invariant metric of
$\SO(n)$ and $g_t(X,Y)=Q(X,P_t(Y))$ for $X,Y\in\fn$. Since
$\h(L)=0$, it follows as in the previous section that $\delta(t)>0$
when $t>0$. Notice that $[Y_2, X_2] -  [X_1, Y_1]= 2
\sum_{s=1}^{s=k-1} a_{s,k-1}E_{s,k}$ if $a=(a_{i,j})\in \SO(k-1)$
and is hence an arbitrary vector in the last column of the upper
$k\times k$-block. Since this column can be regarded as the tangent
space of $\SO(k)/\SO(k-1)=\Sph^{k-1}$ and since by assumption $\gK'$
does not act transitively on $\Sph^{k-1}$, we can  choose
$a\in\SO(k-1)$ such that $[Y_2, X_2] - [X_1, Y_1]$ is perpendicular
to the Lie algebra of $\gK_-$ with respect to $Q$. This implies that
there exists $C>0$ and $\epsilon>0$ such that $Q([Y_2, X_2] -  [X_1,
Y_1],
 P^{-1}_t( [Y_2, X_2] -  [X_1, Y_1] ) )\le C$ for $t\in
(0,\epsilon)$. Hence  \eqref{Bound} implies that $\log(\delta)'\le
C$ which again contradicts $\delta(0)=0$ and $\delta(t)>0$ for
$0<t<\epsilon$.
\end{proof}

 By the classification of transitive actions on spheres, condition
(b) in \tref{newones}  only excludes $\mu=sl$ and  the spin
representations of $\Spin(7)$ and $\Spin(9)$.
 But notice that in the
case of $\mu=sl$, the above cohomogeneity one manifold $M$ does
admit an invariant non-negatively curved metric, since it is
diffeomorphic to the homogeneous space $\SO(n+1)/\mu(\gK)\cdot
\SO(n-k+1)$ endowed with the cohomogeneity one action of
$\SO(n)\subset \SO(n+1)$.

On the other hand, condition (a) is very restrictive. For
$\K'/\H'=\SO(\ell)/\SO(\ell -1)$, such representations  are called
\emph{class one representations} of the symmetric pair
$(\SO(\ell),\SO(\ell -1))$  and it is well known (see e.g.
\cite{Wa}) that they consist of the irreducible representations
$\mu_m$ of $\SO(\ell)$ on the
 homogeneous harmonic polynomials of degree $m$ in $\ell$ variables.
It is now easy to see that for these special representations,
condition (c) is  satisfied as well. In particular we obtain:

\begin{cor}\label{maincor}
For each  $\ell\ge 3$  and $m\ge 2$, the pair
$(\gK',\gH')=(\SO(\ell),\SO(\ell -1))$ and the representation
$\mu_m$ defines a \com as in \eqref{ex} that does not admit an
invariant metric with non-negative curvature.

\end{cor}

In these examples the codimension of the singular orbits are
$\ell\ge 3$ and $n-k+1\ge 3$. In particular, for any pair of
integers $(\ell_-,\ell_+)$ with $\ell_{\pm} \ge 3$, there exists an
infinite family of cohomogeneity one manifold with singular orbits
of codimensions $\ell_{\pm}$ and no invariant metric of non-negative
sectional curvature. This, together with Theorem A,  implies Theorem
B.

For a general  $\gK'/\gH'=$$ \Sph^{\ell -1}$ it is clear that
$S^m(sl)$ restricted to  $\gH'$ has a fixed vector. This implies
that at least one of the irreducible subrepresentations $\mu$ of
$S^m(sl)$ also satisfies condition (a) and one again shows that
condition (c) is satisfied. This gives rise to at least one
representation satisfying the conditions of \tref{newones} for any
$\K'/\H'$.

\providecommand{\bysame}{\leavevmode\hbox
 to3em{\hrulefill}\thinspace}

 \end{document}